\documentclass[12pt]{article}
\usepackage{amsmath}
\usepackage{graphicx,psfrag,epsf}
\usepackage{natbib}
\usepackage{amsmath,amsfonts,amsthm}
\usepackage{hyperref,todonotes}
\usepackage{url} 
\usepackage{tikz}
\usetikzlibrary{shapes.geometric, arrows.meta, positioning, fit}
\usepackage{algorithm}
\usepackage{algpseudocode}
\usepackage{soul}

\usepackage{amscd}
\usepackage{xypic}
\input{xy}
\xyoption{all}
\xyoption{poly}
\usepackage[all]{xy}
\usepackage{tikz}
\newtheorem{theorem}{Theorem}[section]

\newtheorem{remark}{Remark}

\tikzset{
    block/.style = {
        rectangle,
        draw,
        fill=white,
        text width=6em,
        align=center,
        rounded corners,
        minimum height=3em
    },
    line/.style = {
        draw, -Latex
    },
    group/.style = {
        draw,
        dashed,
        rectangle,
        inner sep=2mm
    }
}

\begin{document}

\def\spacingset#1{\renewcommand{\baselinestretch}%
{#1}\small\normalsize} \spacingset{1}


{
  \title{\bf Distinguishing Calabi-Yau Topology using Machine Learning}
  \author{Yang-Hui He
    \hspace{.2cm}\\
    {\small
    London Institute for Mathematical Sciences, Royal Institution, W1S 4BS, UK}\\
    {\small
    Merton College, University of Oxford, OX14JD, UK
    }\\
    Zhi-Gang Yao \\
    {\small
    Department of Statistics and Data Science, National University of Singapore, Singapore
    }\\
    Shing-Tung Yau \\
    {\small
    Yau Mathematical Sciences Center, Tsinghua University, Beijing, China}\\
    ~\\
    {\small
    hey@maths.ox.ac.uk, \
    zhigang.yao@nus.edu.sg, \ 
    styau@tsinghua.edu.cn
    }
    }
  \maketitle

\bigskip
\begin{abstract}
While the earliest applications of AI methodologies to pure mathematics and theoretical physics began with the study of Hodge numbers of Calabi-Yau manifolds, the topology type of such manifold also crucially depend on their intersection theory.
Continuing the paradigm of machine learning algebraic geometry, we here investigate the triple intersection numbers, focusing on certain divisibility invariants constructed therefrom, using the Inception convolutional neural network.
We find $\sim 90\%$ accuracies in prediction in a standard fivefold cross-validation, signifying that more sophisticated tasks of identification of manifold topologies can also be performed by machine learning.
\end{abstract}

\noindent%
\vfill

\newpage
\spacingset{1.45} 

\section{Introduction}
In the realm of contemporary scientific inquiry, the application of deep learning has emerged as a transformative force, endowing researchers with novel tools and methodologies to explore intricate scientific phenomena. 
While AI has been transformative in the experimental sciences for decades, its revolution in theoretical and mathematical sciences have been a relatively more recent emergence \citep{He:2017aed,Krefl:2017yox,Carifio:2017bov,Ruehle:2017mzq}; see reviews and prospects in \citet{He:2018jtw,Ruehle:2020jrk,He:2021oav,Bao:2022rup,he2024can,gukov2024rigor,He:2024gnk}.
A programme to machine-learn various structures in mathematics by looking at {\it pure mathematical data} has been launched over the past seven years, and has blossomed to include algebraic geometry \citep{He:2017aed,He:2018jtw}, algebra and representation theory \citep{He:2019nzx,davies2021advancing}, graph theory and combinatorics \citep{He:2020fdg}, knot theory \citep{Craven:2020bdz,Gukov:2020qaj,davies2021advancing}, number theory \citep{Alessandretti:2019jbs,He:2022pqn}, symbolic computation \citep{peifer2020learning,England_2018,lample2019deep}, etc.

Calabi-Yau manifolds are a protagonist in modern mathematics, residing at the interface between algebraic and differential geometry, mathematical physics (especially string theory), topology, number theory, and dynamical systems.
These are K\"ahler manifolds that admit zero Ricci curvature.
The simplest example is one that is familiar to the beginning student: the torus $T^2 = S^1 \times S^1$ can be complexified into a Riemann surface, which is in turn trivially K\"ahler and furthermore is Ricci-flat.
The name ``Calabi-Yau'' originated in the Fields-winning proof by one of the authors (STY) \citet{yau1977calabi} that settled the Calabi Conjecture \citep{calabi1957kahler}, and was first coined by physicists \citet{candelas1985vacuum}.

In searching for the Standard Model of particle physics within string theory compactifications, one of the first datasets in modern geometry was created; this was the so-called CICYs, for ``Complete Intersection Calabi-Yau threefolds'' \citep{candelas1988complete,gagnon1994exhaustive,green1989all} (for pedagogical introductions, see the classic in \citet{hubsch1992calabi}, as well as a recent treatment \citet{He:2018jtw}).
These are Calabi-Yau manifolds of complex dimension three, realized as algebraic varieties embedded as complete intersections in products of projective spaces.
From a data-scientific point of view, they are rather straight-forwardly represented as matrices of non-negative integers, which encode the multi-degree information of the homogeneous polynomials that define the variety.
It was natural, therefore, that the first machine-learning experiments in pure mathematics was in consideration on this data-set \citep{He:2017aed}.
Indeed, despite its rich properties, the precise nature of the Calabi-Yau manifold often proves challenging to intuitively grasp and compute, presenting an intriguing and demanding avenue for the application of deep learning. 

Now, one should bear in mind that in topology, characteristics that describe the structure of manifolds, such as the classic Betti numbers which describe the closed cycles, are often intricately and closely linked with the data representing these manifolds. This makes the precise derivation of these characteristics quite challenging. The analytic mathematical formulas for these properties are largely unexplored \footnote{Due to the moduli dependence of quantities such as cohomology groups—unlike their alternating sums, which give rise to Euler characteristics that can be more readily computed by the index theorem—they are very difficult to calculate. See \citet{Constantin:2018hvl} for a recent successful extraction of exact formulae for surfaces guided by machine learning.}, and their computation relies on algorithms that are both complex and time-consuming. Therefore, the application of deep learning in the field of topology poses significant challenges. This paper is dedicated to utilizing the capabilities of deep neural networks to predict some key properties of the Calabi-Yau manifold, in order to explore the potential applications of deep learning within the realm of mathematics and theoretical science.


Henceforth, we focus on CICY threefolds, which have been a pivotal class of manifolds in the context of string model building and, as mentioned, one of the first data-set in pure mathematics to enter the age of data \citep{candelas1988complete} and the era of machine-learning \citep{He:2017aed}.
Since then, there has been a host of activity on machine-learning properties of CICYs.
Neural network methodologies were applied to predict the Hodge numbers of CICY three-folds, favorability, and discrete symmetries \citep{bull2018machine,Bull:2019cij}, indicating that they outperform classical statistical methods such as Support Vector Machines (SVM).
Subsequently, the prediction accuracy of Hodge numbers were increased to over 99\% by leveraging the Inception network model  \citep{erbin2021inception,Erbin:2020tks}, and more recent comparative study, with impressive accuracies, was performed in \citet{Keita:2024skh}.
Other related directions have included
CICY fourfolds \citep{he2021machine,Erbin:2021hmx} 
and generalized CICYs \citep{Cui:2022cxe}
(see review in \citet{Erbin:2023ncy}), numerical metrics \citep{Ashmore:2019wzb,Douglas:2020hpv,Anderson:2020hux,Anderson:2023viv,Ashmore:2023ajy}, weighted CICYs \citep{Berman:2021mcw,Hirst:2023kdl,macfadden2024dna}, distinguishing elliptic fibrations \citep{He:2019vsj} and topologies \citep{Jejjala:2022lxh, chandra2023manifolds, gendler2023counting}, etc.
This expansion reflects a continued effort to broaden the scope of neural network applications in predicting geometric properties, highlighting the versatility of these models in handling increasingly complex mathematical structures.

Now, most of the work on machine-learning of CICY topological invariants have focused on Hodge numbers, whereas, as we will soon see, more refined invariants involve Chern classes and intersection numbers.
This is an important issue because the ``topological type'' of a (smooth, compact, and simply connected) Calabi-Yau threefold, due to the extension of a theorem of Wall \citet{wall1966classification}, consists of a pair of Hodge numbers $(h^{1,1}, h^{2,1})$, as well as (when expressed in an appropriate basis of the K\"ahler classes) the second Chern class $(c_2)_r$ and the triple-intersection numbers $d_{rst}$ of the curve classes. 
It is a standing conjecture of one of the authors (STY), that in any complex dimension $n$, the possible topological types of a Calabi-Yau $n$-fold is {\it finite}.

In sum, we hope to explore and expand the potential of deep learning to predict more critical properties of the CICY three-folds. In this study, we seek to predict the triple intersection numbers of CICY three-folds, combining ideas from \citet{He:2017aed,erbin2021inception,Jejjala:2022lxh}, by an {\it inception network V3 model}. By amalgamating the nonlinear feature extraction capabilities of deep learning with the geometric properties of manifolds, we aim to enhance the deep learning model's understanding of the CICY three-folds.

\subsection*{Summary}
The organization of the paper is as follows.
We begin in Section 2 with a review of the construction of the complete intersection Calabi-Yau (CICY) database, emphasizing on their topological invariants, especially the triple-intersection numbers.
Then, in Section 3, we compare four models - SVM, random forest, XGBoost, and Google Inception, in predicting the four key quantities $(d_1, d_2, d_3, d_p)$ defined in \eqref{d-def} which capture CICY topology.
In Section 4, we report the prediction accuracies of these quantities, with the main conclusion being that the Inception Network can reach around 90\% accuracy in the standard 80-20 data cross-split.
Finally, we conclude with outlook in Section 5.

\subsection*{Acknowledgments}
YHH is grateful to the STFC, UK for grant
ST/J00037X/2 and a Leverhulme project grant.
ZGY has been supported by Singapore Ministry of Education Tier 2 grant A-8001562-00-00 and Tier 1 grant A-8002931-00-00 at the National University of Singapore. ZGY also thanks the London Mathematical Society and the Royal Statistical Society for their partial support of his research trip to the UK.

\section{Background and Dataset}
We begin with a brief description of the construction of CICY threefolds, as well as their topological properties.
Throughout, we will focus on the representation of the data.

\subsection{Dataset}
CICY manifolds are realized as a complete intersection of polynomials in a product of complex projective spaces \citep{candelas1988complete,gagnon1994exhaustive}. 
What this means is the following.
Consider the ambient space as the product $X=\mathbb{C}P^{n_1}\times\cdots\times\mathbb{C}P^{n_m}$, a CICY three-fold $M^3$ is embedded as $k$ polynomials, the complete intersection condition is that
\begin{equation}
    \sum_{i=1}^m n_i = 3 + k 
\end{equation}
so that $\dim(X) - \dim(M^3) = 3$.
We will index the projective space factors by $i = 1, 2, \ldots, m$ and the polynomials by $r = 1, 2, \ldots, k$.
The $r$-th defining polynomial is then homogeneous of multi-degree $q_i^r$ with respect to $\mathbb{C}P^{n_i}$; this can be recorded as a configuration matrix
\begin{equation}\label{amatrix}
M^3 \simeq 
A_{m\times k}=\begin{pmatrix}
q_1^1& \cdots&q_1^k\\
\vdots& \ddots&\vdots\\
q_m^1& \cdots&q_m^k
\end{pmatrix}
\end{equation}
Thus, each deformation family of CICYs is represented by an $m \times k$ matrix of non-negative (and indeed largely sparse) integers.

Now, $M^3$ is a Calabi-Yau three-fold, so it has vanishing first Chern class. 
This conveniently translates to (using adjunction, see e.g., \citet{hubsch1992calabi,He:2018jtw}) the condition that row sums to one less than dimension of the corresponding projective factor:
\begin{equation}
A\cdot \begin{pmatrix}
1\\
\vdots\\
1
\end{pmatrix}_{k\times 1}=\begin{pmatrix}
n_1+1\\
\vdots\\
n_m+1
\end{pmatrix}.    
\end{equation}

The CICY three-folds were classified in \citet{candelas1988complete}, using the then state-of-the-art computing.
Up to trivial permutation equivalence and also simple equivalence due to birational transformations, there are 7890 configuration matrices, giving rise to perhaps the first ``big'' database in geometry.
Today, this database was updated in \citet{Anderson:2007nc} and maintained at
\begin{quote}
\url{http://www-thphys.physics.ox.ac.uk/projects/CalabiYau/cicylist/}
\end{quote}

Now, an integral cohomology ring 
\begin{equation}
H^*(\mathbb{C}P^{n_1}\times\cdots\times\mathbb{C}P^{n_m};\mathbb{Z} )\cong \mathbb{Z}[x_1,\cdots,x_m]/\langle x_1^{n_1+1},\cdots, x_m^{n_m+1} \rangle    
\end{equation}
descends from the ambient $X$, and
and we have, correspondingly, $k$ cohomology classes 
$\sum\limits_{i=1}^m q_i^j x_i$ for $j = 1, \ldots, k$.
It should be emphasized that this implies that while $h^{1,1}(M^3)$ of these classes come at least from the $m$ projective factors, by no means they are equal in general, as extra K\"ahler classes arise from the restriction.
In the case of equality $h^{1,1}(M^3) = m$, when {\it all} divisors of $M^3$ descend from the simple ambient space $X$, we refer to $M^3$ as {\it favourable}. Of the original 7890 CICY threefold configuration matrices \citep{candelas1988complete}, there are 4896 favorable geometries and 2994 unfavorable geometries \citep{anderson2017fibrations}.  Furthermore, new and favourable descriptions of 2946  unfavorable CICYs were presented therein. The favorable list can be obtained from \url{http://www1.phys.vt.edu/cicydata/}.

Here, we briefly summarize $(d_1,d_2,d_3,d_p)$ for all 7,680 samples with the empirical distributions:

\begin{itemize}
    \item \textbf{$d_1$:} takes values in $\{1,2,3,4,5,6,8,9,12,16\}$. The vast majority are $d_1=1$ (7336 samples), with a smaller cluster at $d_1=2$ (413). Other values are rare: e.g.\ $d_1=3$ (24), $d_1=4$ (34), and only isolated cases at $d_1=5,6,8,9,12,16$.
    \item \textbf{$d_2$:} shares the same range $\{1,2,3,4,5,6,8,9,12,16\}$, but exhibits a different balance: $d_2=1$ (5315) and $d_2=2$ (2402) dominate, while higher values (3–16) occur with very low frequency.
    \item \textbf{$d_3$:} has extended range $\{1,2,3,4,5,6,8,9,12,16,18\}$. Its distribution is sharply peaked at $d_3=2$ (3728) and $d_3=6$ (3710), almost evenly split. All other values appear in very small counts.
    \item \textbf{$d_p$:} shows the broadest spread, covering $\{2,4,6,8,12,18,24,36,44,50,52,54,56,60,64\}$. It is, however, highly concentrated at $d_p=4$ (3724) and $d_p=12$ (3688). Smaller groups include $d_p=2$ (275), $d_p=6$ (25), $d_p=8$ (61), and a tail of rare larger values up to $64$.
\end{itemize}

Now, a classic theorem in K\"ahler geometry is:
\begin{theorem}[C.~T.~C.~Wall \citep{wall1966classification}]
The topological type of a compact K\"ahler threefold is completely determined by 
\begin{enumerate}
	\item the Hodge numbers $h^{p,q}$;
	\item the triple intersection numbers $d_{rst}$;
	\item the first Pontrjagin class $p_1=c_1^2-2c_2$.
\end{enumerate}
\end{theorem}
For Calabi-Yau threefolds, $c_1 = 0$ by definition.
Furthermore, for the second Chern class, we can fix a K\"ahler basis $\{ J^r \}_{r = 1,2,\ldots,h^{1,1}}$ for $H^2(M^3; \mathbb{Z})$, into which one can expand $c_2 = \sum\limits_{s,t = 1}^{h^{1,1}}(c_2)_{st} J^s J^t$.
Likewise, the triple intersection numbers record the information about $H^1(M^3; \mathbb{Z}) \times H^1(M^3; \mathbb{Z}) \times H^1(M^3; \mathbb{Z}) \longrightarrow \mathbb{C}$, can be expressed in this basis as the triple integral $d_{rst} = \int_{M^3} J^r \wedge J^s \wedge J^t$.
Moreover, for simply-connected Calabi-Yau threefolds (with which we focus here, and to which class all CICY certainly belong), the only non-trivial Hodge numbers are $(h^{1,1}, h^{2,1})$ (these were computed in \citet{green1989all}).
Therefore, the topological type of Calabi-Yau threefolds will be given by the following list of non-negative integers:
\begin{equation}
 Y:=
 \{
 h^{1,1}, h^{2,1}, d_{rst}, (c_2)_{st}  
 \}   
\end{equation}
This list of integer (tensors) is the topological datum for our Calabi-Yau threefold.

 

\subsection{CICY Intersections}
 As mentioned above, for a favorable CICY, one identifies the K\"ahler basis 
 $J^r$ as coming from the ambient projective space and thus set
 $\{J^r = x^r\}_{r=1,\cdots, h^{1,1}}$ for $H^2(M^3,\mathbb{Z})$. The triple intersection form in this basis is written as 
\begin{align}
\nonumber
d_{rst}
& := \int_{M} x^r\wedge x^s \wedge x^t 
 =\langle x_rx_sx_t, [M^3]  \rangle \\
& =\langle x_rx_sx_t\cdot \prod_{j=1}^{k}(q_1^jx_1+\cdots+q_m^jx_m) , [\mathbb{C}P^{n_1}\times\cdots\times\mathbb{C}P^{n_m}] \rangle
\end{align}
where in the second line one computes it explicitly by pulling back integration from the ambient space $X$.

Similarly, the Chern class of $M^3$ is
\begin{align}
\nonumber
c(M) & = \prod_{i=1}^{m}(1+x_i)^{n_i+1}\cdot\prod_{j=1}^{k}(1+q_1^jx_1+\cdots+q_m^jx_m)^{-1} \\
& =1+\sum_{r,s}[c_2(M)]_{rs}x_rx_s+\sum_{r,s,t}[c_3(M)]_{rst}x_rx_sx_t \ ,
\end{align}
when expanded into our basis, for $r,s,t=1,\cdots, h^{1,1}$. 
Explicitly, we can express them in terms of the configuration matrix $A$:
\begin{align}
\nonumber
c_1(M^3) & = 0 \\
\nonumber
[c_2(M^3)]_{rs}
&=\frac{1}{2}\left[-\delta_{rs}(n_r+1)+\sum_{j=1}^{h_{11}}q_j^rq_j^s  \right] \\
[c_3(M^3)]_{rst} &=
\frac{1}{3}\left[\delta_{rst}(n_r+1)-\sum_{j=1}^{h_{11}}q_j^rq_j^sq_j^t  \right]  \ .
\end{align}
Now, the Euler number is easily determined by the intersection numbers as
\begin{equation}
    \chi(M^3) = \sum_{r,s,t=1}^{h_{11}} d_{rst} [c_3]_{rst} = 2(h^{1,1} - h^{2,1}) \ ,    
\end{equation}
where the first equality follows from Gauss-Bonnet-Chern in expressing the Euler number in terms of the top Chern class, and the second equality come from expressing the Euler number in terms of an alternating sum over the Hodge (Betti) numbers.
Also, by Poinca\'re duality, $H^4(M,\mathbb{Z})\cong H^2(M,\mathbb{Z})$, thus second Chern class $c_2$ can be labelled by $h^{1,1}$ as 
\begin{equation}
[c_2(M)]_t=\sum_{r,s}[c_2(M)]_{rs}d_{rst}.    
\end{equation}

Therefore, there are two parts to the problem when distinguishing CICY manifolds, the first is to compute the data $(d_{rst}, c_2)$ for each manifold in the list, the second is to compare the resulting quantities and to decide when they correspond to different manifolds. 

In general the computation of $(d_{rst}, c_2)$ is applicable for favorable geometries, we will design the computer algorithm later.
On the other hand, it is usually hard to determine whether two sets of data $(d_{rst}, c_2)$ corresponding to the same topology, i.e., to  determine whether the difference between them corresponding merely to a change of basis of $H^2(M,\mathbb{Z})$. 
Thus, we adopt a family of divisibilitiy invariants (see \S 8.1 of \citet{hubsch1992calabi}) up to change of basis:
\begin{align}
\nonumber
d_1 & := \gcd  \{d_{rst}\} ; \\
\nonumber
d_2 & := \gcd  \{d_{rrs}, 2d_{rst}\}; \\
\nonumber
d_3 & := \gcd \{d_{rrr}, 3 (d_{rrs}\pm d_{rss}), 6 d_{rst}\}; \\
\label{d-def}
d_p & :=\gcd \{[c_2]_{r}\}.
\end{align}
where $\gcd$ is the greatest common divisor.
In this way, the new indices 
\begin{equation}
Y':=\{h^{1,1}, h^{2,1}, d_1,d_2,d_3, d_p  \}  
\end{equation}
will be used to distinguish CICY topology. 
Since Hodge number prediction had been the subject of the literature since 2017, we will here focus on $(d_1,d_2,d_3,d_p)$.
To our knowledge, the compilation of these four quantities for manifolds is a new addition to existing databases.

\subsection{Computer Algorithm for $d_{rst}$}
As promised, let us digress briefly to present an algorithm for computing $d_{rst}$ for favorable CICY's.
Recall that the CICY configuration matrix $A = (q_m^j)$ defined in \eqref{amatrix} satisfies
\begin{equation}
  \sum_{i=1}^m n_i =3+k;  \ k\leq 18, \ m\leq 15 \ ,   \qquad
  {A}\cdot \begin{bmatrix}
1\\
\vdots\\
1
\end{bmatrix}_{k\times 1}=\begin{bmatrix}
n_1+1\\
\vdots\\
n_m+1
\end{bmatrix}.
\end{equation}
where the $A$ can also be obtained from the data in \url{http://www1.phys.vt.edu/cicydata/}.

We have the explicit algorithm as follows
\begin{enumerate}
\item
First, we can write the intersection form in terms of the configuration matrix entries (the multi-degrees that define the complete intersection) as
\begin{align}
\nonumber
d_{rst}&=\langle x_rx_sx_t, [M]  \rangle=\langle x_rx_sx_t\cdot \prod_{j=1}^{k}(q_1^jx_1+\cdots+q_m^jx_m) , [\mathbb{C}P^{n_1}\times\cdots\times\mathbb{C}P^{n_m}] \rangle\\
\nonumber
&=  \text{coefficient of}\ \prod_{i=1}^{m}x_i^{n_i} \ \text{in}  \ x_rx_sx_t\cdot \prod_{j=1}^{k}(q_1^jx_1+\cdots+q_m^jx_m)\\
&=\sum_{\sigma} q_1^{\sigma_1}q_2^{\sigma_2}\cdots q_m^{\sigma_m} \ ,
\end{align}
where $\sigma=(\sigma_1,\sigma_2,\cdots,\sigma_m)$ is the partition of set $[k]=\{1,\cdots, k \} $, $q_i^{\sigma_i}:=\prod_{j\in \sigma_i} q_i^j$. $|\sigma_i|=n_i,\ i\neq r, s, t$; $|\sigma_i|=n_i-1 $ otherwise.  

\item

{Second, we construct a $k\times k$ square matrix $A_{rst}(i,j)$  by adding $n_i-1$ copies of the $i$-th row to $A$ if $i\neq r, s, t$. Otherwise, we add  $n_i-2$ copies of the $i$-th row to $A$}. 
Note that the subscript $rst$ is just to show that the matrix depends on the indices of the triple intersection numbers $d_{rst}$ and the row-column indices that actually define the matrix are denoted as $(i,j)$.

In this way, we can more conveniently write $d_{rst}$ essentially as a weighted permanent:
\begin{equation}
d_{rst}=\sum_{\sigma\in S_k} \prod_{i=1}^k A_{rst}(i, \sigma(i))/ (\prod_{j\neq r,s,t} n_j !\cdot \prod_{j= r,s,t} (n_j-1) !) \ .    
\end{equation}
Note that we can also define 
\begin{equation}
\text{detm}(A_{rst})=\sum_{\sigma\in S_k} \prod_{i=1}^k A_{rst}(i, \sigma(i)) \ ,    
\end{equation}
which is the permanent, the unsigned version of the determinant (the code can be modified therefrom). 

{
\begin{remark}
Note that in special case when $n_i=1$, ``adding  $n_i-2$ copies of the $i$-th row" means ``deleting the $i$-th row". Thus it is possible that the size of $A_{rst}$ is smaller than $A$.
\end{remark}}

The pseudo code is shown in Algorithm 1. 

\item
Finally, store these $d_{rst}$ and calculate $c_2$.
\end{enumerate}

\begin{algorithm}
\caption{An algorithm to calculate detm(A)}
\begin{algorithmic}[1]
\Function{detm}{$A$}
    \If {number of rows in A $= 1$}
        \State $f \gets A(1, 1)$
    \ElsIf {number of rows in A $= 2$}
        \State $f \gets A(1, 1) \times A(2, 2) - A(1, 2) \times A(2, 1)$
    \Else
        \State $f \gets 0$
        \For {$i \gets 1$ \textbf{to} number of rows in A}
            \State $B \gets A$
            \State \Call{Remove Row}{$B, i$}
            \State \Call{Remove Column}{$B, 1$}
            \If {$A(i, 1) \neq 0$}
                \State $f \gets f + A(i, 1) \times \Call{detm}{B}$
            \EndIf
        \EndFor
    \EndIf
    \State \Return $f$
\EndFunction
\end{algorithmic}
\end{algorithm}


\section{Methods}
\subsection{Inception Network model} 
    To construct a robust architecture, we leverage the Inception Neural Network Model from GoogleNet \citep{szegedy2016rethinking}. This deep convolutional neural network architecture employs inception modules, a unique design that captures features across various spatial scales and complexities concurrently, departing from the sequential approach of traditional Convolutional Neural Networks. The distinctive feature of Inception is its ability to efficiently gather both local and global contextual information, resulting in heightened accuracy and diminished computational complexity. Notably, the Inception Network model relieves users of the burden of manually deciding on kernel selection or the necessity of pooling layers. The network autonomously determines these parameters. Users can input a comprehensive set of potential values for these parameters into the network, allowing it to autonomously learn the optimal parameters and convolution kernel combinations. Consequently, inception modules can be stacked iteratively to construct a more extensive network, effectively expanding both its depth and width.

\subsection{Architectures}     
   Our model is displayed in the Figure \ref{fig:Model Structure}, which takes the configuration matrix of the CICY manifolds as input and outputs a 1$\times$4 vector, representing the predicted CICY triple intersection numbers $d_1$, $d_2$, $d_3$ and $d_p$. The architecture consists of three inception modules with 32, 64, and 32 filters, respectively. Within each inception module, we employ the rows 12$\times$1 and the columns 1$\times$15 kernels for two parallel convolutions with a stride of 1 and concatenate the outputs of both layers together over the channel dimension. The Leaky ReLU function with negative slope 0.2 serves as the activation function. 
   
   To maintain consistent data dimensions, we attempt zero-padding and adopt the same padding options. Batch normalization layers follow each concatenation layer, facilitating the connection between these inception modules. The output from the final inception module is passed through a flatten layer, which then connects to a fully connected layer that outputs a four-dimensional vector. During training, we utilize the Adam gradient descent optimizer with Mean Square Error (MSE) as the loss function. The initial learning rate is set at 0.0001, and the batch size was specified as 32. To mitigate overfitting, we introduce $l_2$ regularization with weight 0.002. This comprehensive approach aims to enhance the model's predictive accuracy while addressing potential overfitting challenges. The model was trained over 3000 epochs to ensure thorough learning and stable optimization of the parameters, during which we observed consistent improvements in performance metrics until convergence.

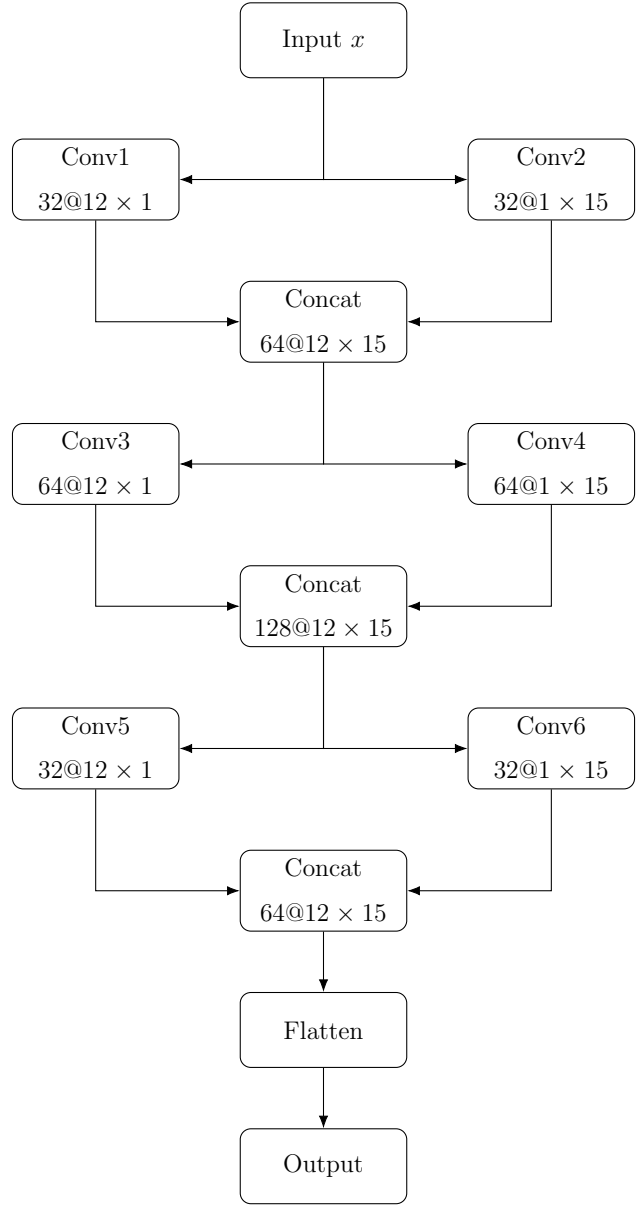
\begin{figure}[h!]
\centering
\begin{tikzpicture}[node distance=1cm and 1cm, auto, scale=0.8, transform shape]
    \node (input) [block] {Input $x$};
    \node (conv1) [block, below left=of input] {Conv1\\$32@12 \times 1$};
    \node (conv2) [block, below right=of input] {Conv2\\$32@1 \times 15$};
    \node (concat1) [block, below right=of conv1] {Concat\\$64@12 \times 15$};
    \node (conv3) [block, below left=of concat1] {Conv3\\$64@12 \times 1$};
    \node (conv4) [block, below right=of concat1] {Conv4\\$64@1 \times 15$};
    \node (concat2) [block, below right=of conv3] {Concat\\$128@12 \times 15$};
    \node (conv5) [block, below left=of concat2] {Conv5\\$32@12 \times 1$};
    \node (conv6) [block, below right=of concat2] {Conv6\\$32@1 \times 15$};
    \node (concat3) [block, below right=of conv5] {Concat\\$64@12 \times 15$};
    \node (flatten) [block, below=of concat3] {Flatten};
    \node (output) [block, below=of flatten] {Output};

    \path[line] (input) |- (conv1);
    \path[line] (input) |- (conv2);
    \path[line] (conv1) |- (concat1);
    \path[line] (conv2) |- (concat1);
    \path[line] (concat1) |- (conv3);
    \path[line] (concat1) |- (conv4);
    \path[line] (conv3) |- (concat2);
    \path[line] (conv4) |- (concat2);
    \path[line] (concat2) |- (conv5);
    \path[line] (concat2) |- (conv6);
    \path[line] (conv5) |- (concat3);
    \path[line] (conv6) |- (concat3);
    \path[line] (concat3) -- (flatten);
    \path[line] (flatten) -- (output);

\end{tikzpicture}
\caption{Model Structure}
\label{fig:Model Structure}
\end{figure}

\section{Results}

    \begin{figure}[htbp]
    	\centering
    	\includegraphics[width = 0.8\textwidth]{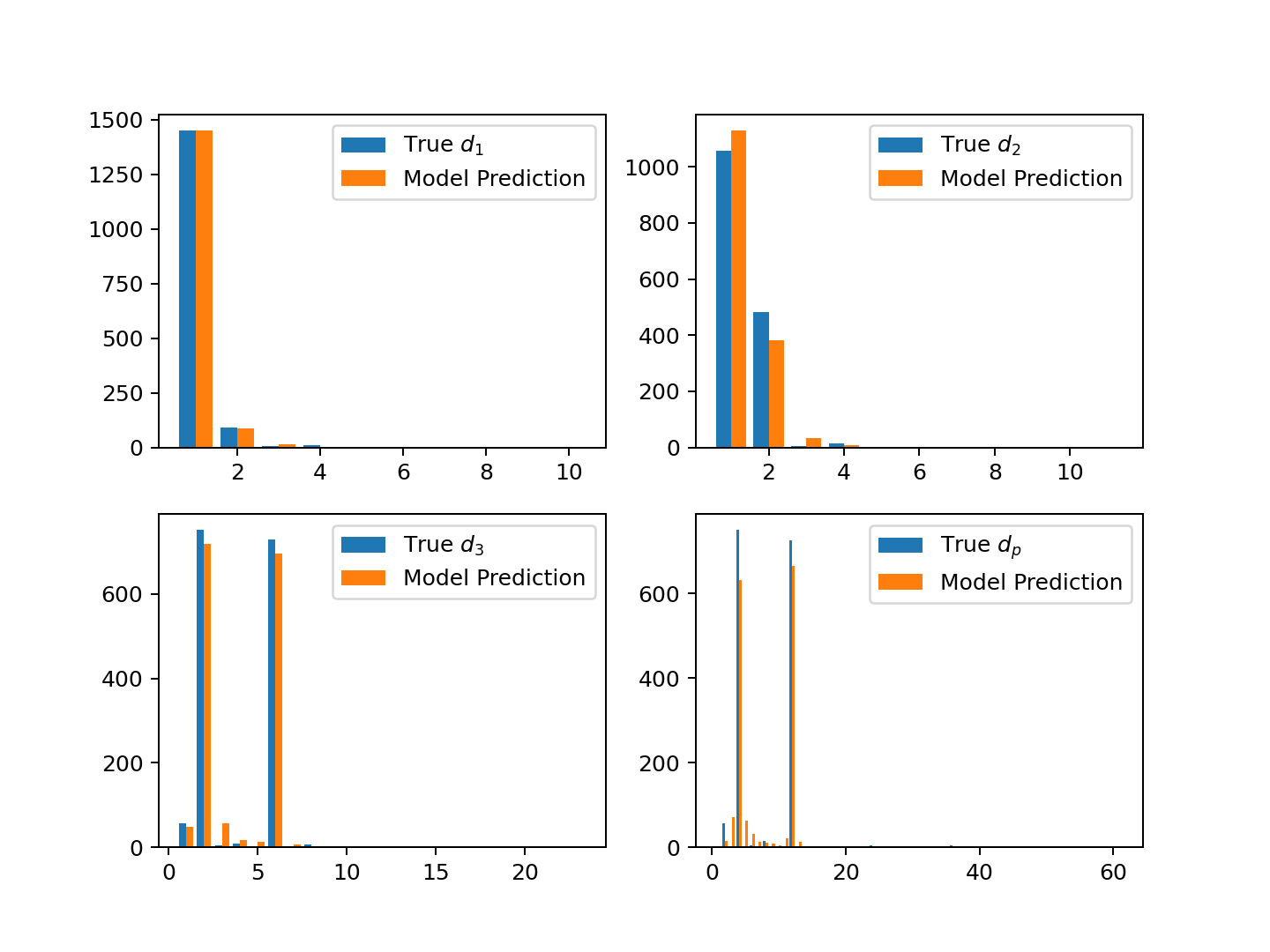}
    	\caption{This depicts comparisons between actual values and model's predictions across CICY triple intersection numbers $d_1$, $d_2$, $d_3$ and $d_p$. We use 80\% of the data as the training set and 20\% of the data as the testing set.}
    	\label{fig:Model Prediction}
    \end{figure}

    \begin{figure}[htbp]
    	\centering
    	\includegraphics[width =0.8\textwidth]{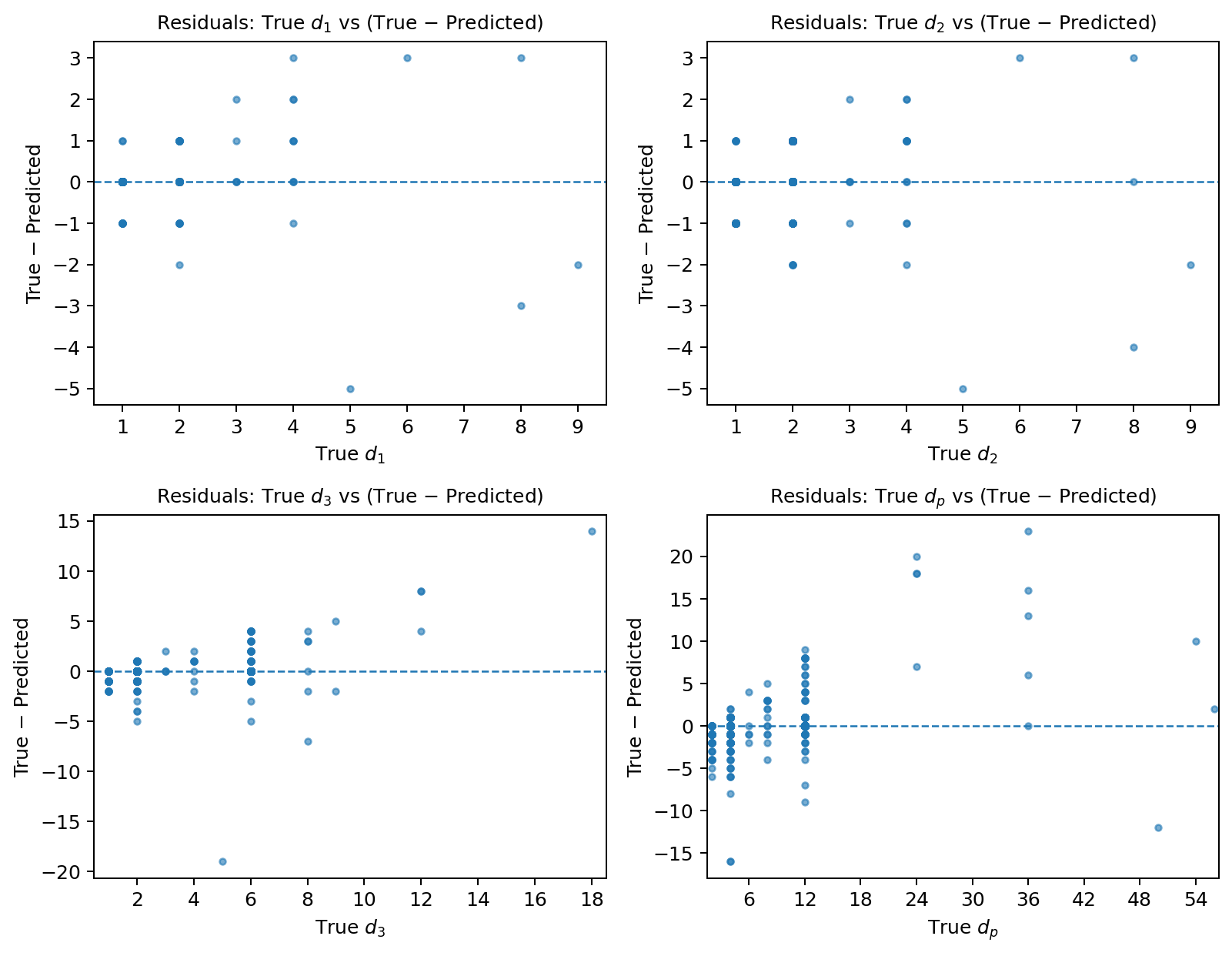}
    	\caption{This residual scatter plots comparing true and predicted CICY triple intersection numbers $d_1$, $d_2$, $d_3$, and $d_p$, based on an 80\%/20\% train–test split.}
    	\label{fig:Model Residuals}
    \end{figure}

    \begin{table}[htbp]
     \centering
     \begin{tabular}{c|ccccc}
  \hline
        {Indicators}   & All  & $d_1$  & $d_2$  & $d_3$ & $d_p$  \\ \hline
   
        {Accuracy}     & 0.886  & 0.971 & 0.838 & 0.907 & 0.827 \\ \
        {MSE}          & 1.039  & 0.074 & 0.215 & 0.808 & 3.060 \\ \hline
     \end{tabular}
    \caption{The accuracy and MSE of our model in predicting the CICY triple intersection numbers on the test set, including overall indicators and indicators for each CICY intersection:  80\% of the data are used as the training set and 20\% of the data as the testing set.}
    \label{tab:Result_1}
   \end{table}

    \begin{figure}[htbp]
    	\centering
              \includegraphics[width=0.45\textwidth]{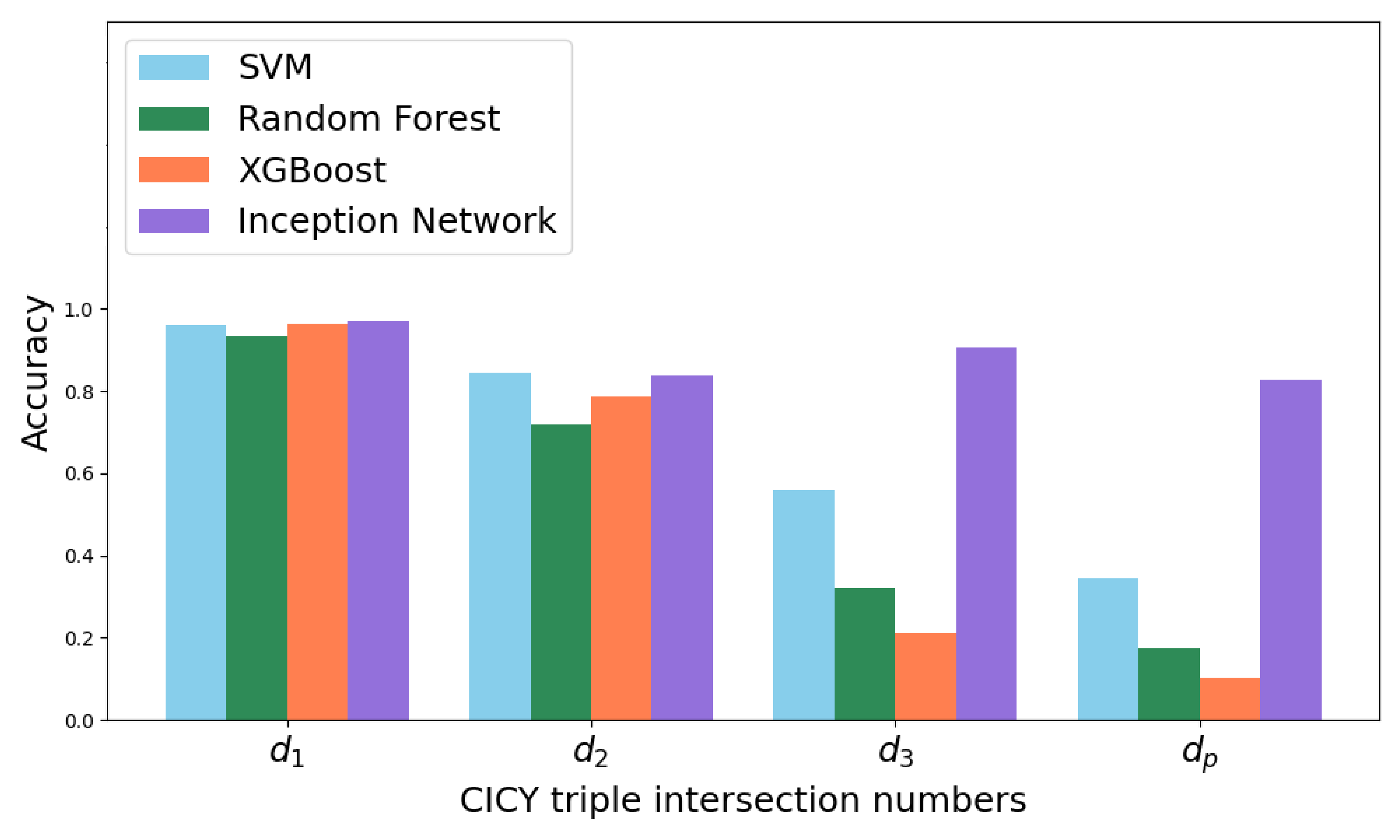}
              \hspace{10pt}
              \includegraphics[width=0.45\textwidth]{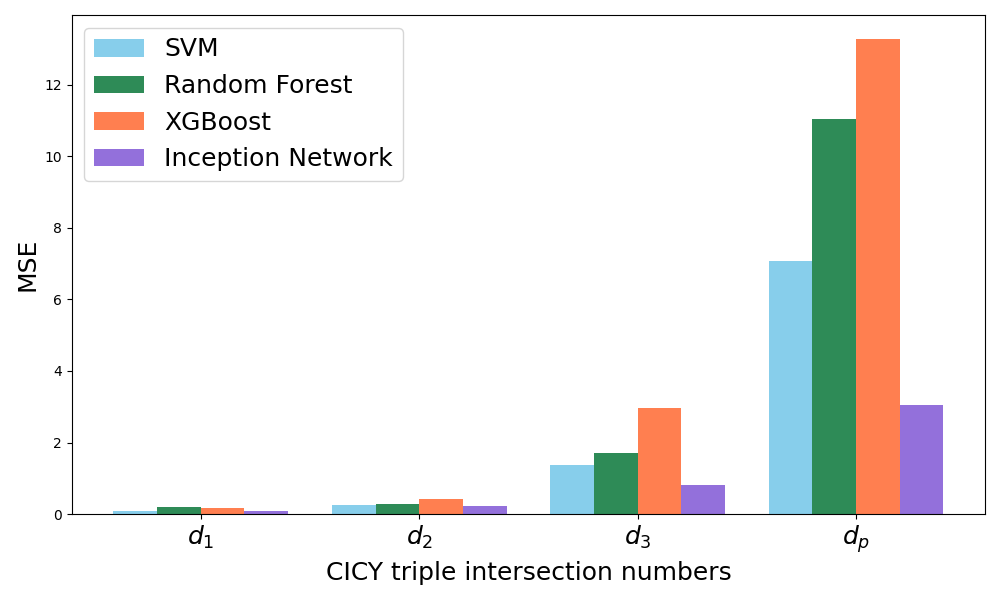}
    	\caption{The accuracy and MSE of different regressors across CICY triple intersection numbers $d_1$, $d_2$, $d_3$ and $d_p$, including SVM regressor, Random Forest regressor, XGBoost regressor and our Inception Network model: 80\% of the data are used as the training set and 20\% of the data as the testing set.}
    	\label{fig:different methods}
    \end{figure}

    \begin{figure}[htbp]
    	\centering
              \includegraphics[width=0.45\linewidth]{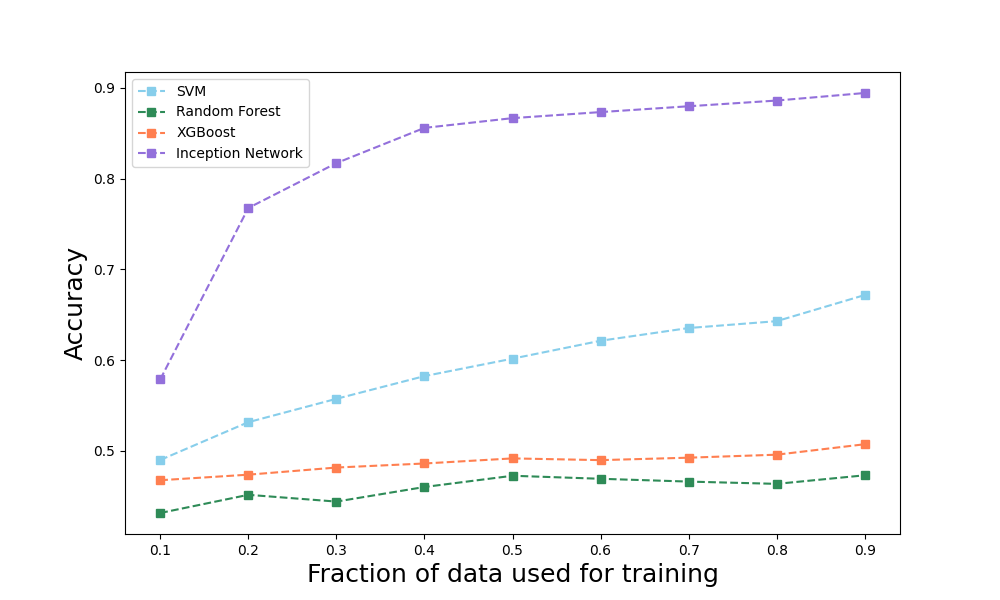}
              \hspace{10pt}
              \includegraphics[width=0.45\linewidth]{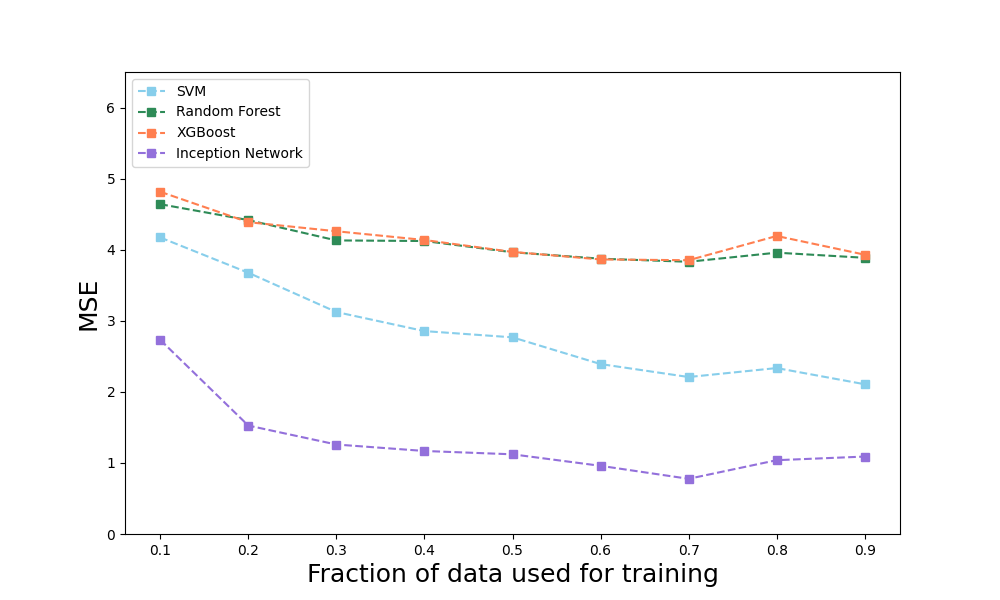}
    	\caption{This presents a comparative analysis between our Inception Network model and conventional machine learning regressors across various fractions of training data. The first figure presents a line graph measuring the accuracy of different regressors in different fractions of training data while the second figure presents MSE.}
    	\label{fig:different methods for different radio of training data}
    \end{figure}

    As proposed in Section 2, we attempt to use machine learning methods to predict CICY triple intersection numbers of  including $d_1$, $d_2$, $d_3$ and $d_p$. The main results demonstrate that our model can perform well in predicting CICY triple intersection numbers, as summarized in Figure
    \ref{fig:Model Prediction} and Table \ref{tab:Result_1}. 
    
    Moreover, we conduct a comprehensive comparison between our deep learning model and several classic regression models. This comparison is aimed at evaluating the relative effectiveness and accuracy of our model in contrast to traditional methods. The results demonstrate that deep learning methods significantly out-perform other traditional regression methods as shown in Figure \ref{fig:different methods}. Furthermore, we delve into an extensive investigation of our model's performance across various ratios of training set splits. This exploration provides insights into the adaptability and robustness of our model under different training conditions. We find that our model is capable of making good predictions across various training set splits compared with other regressors and display the results in Figure \ref{fig:different methods for different radio of training data}.

     Figure \ref{fig:Model Prediction} displays a set of histograms representing the model's predictions versus the true data for the CICY triple intersection numbers $d_1$, $d_2$, $d_3$ and $d_p$. We find that there is a close match between the true data and the model's predictions, with both sets of bars aligning closely, indicating a strong predictive performance, especially in the $d_1$. 

     Figure \ref{fig:Model Residuals} presents residual scatter plots comparing the true and predicted values of the CICY triple intersection numbers $d_1$, $d_2$, $d_3$, and $d_p$. Each subplot displays the deviation $(\text{True} - \text{Predicted})$ with respect to the true values, providing a sample-wise view of model accuracy. The residuals are generally centered around zero, indicating that the model captures the underlying distribution of the data well. Slightly larger deviations are observed for $d_3$ and $d_p$, suggesting that these quantities exhibit higher structural complexity and thus pose greater learning difficulty for the network.

     Table \ref{tab:Result_1} provides a comprehensive analysis of our model's performance on the individual intersection numbers $d_1$, $d_2$, $d_3$, and $d_p$. It highlights that our model achieves high accuracy for each of these CICY triple intersection numbers. Additionally, the MSE across these dimensions is kept low, indicating that the model not only predicts with high accuracy but also with minimal error, thus demonstrating robust performance even in the handling of complex manifold structures. This effective combination of high accuracy and low MSE underscores the efficacy of our model in dealing with the intricacies of CICY manifold predictions. We can observe that $d_p$ exhibits relatively a bit lower prediction accuracy and higher MSE, which stems from the more intricate structure of $d_p$ possibly. In particular, $d_1$ is almost always $1$ or $2$, making its distribution highly concentrated and relatively easy to predict. By contrast, $d_p$ is derived from contractions of the second Chern class and spans a much broader range of values, up to $64$.
    
    In Figure \ref{fig:different methods}, we employed some distinct classic machine learning regression models to predict the CICY triple intersection numbers $d_1$, $d_2$, $d_3$ and $d_p$, and compared their performance with our model. These algorithms include  SVM regressor \citep{cortes1995support}, Random Forest regressor \citep{breiman2001random} and XGBoost regressor \citep{chen2016xgboost}, which served as benchmark models. The settings of these models are as follows.
    
    \textbf{SVM:} The SVM model utilizes the Radial Basis Function (RBF) kernel, a common choice for non-linear data. The regularization parameter, \( C \), is set to 100, which determines the trade-off between achieving a low training error and a low testing error by controlling the model's complexity. The kernel coefficient, \( \gamma \), is set to 0.1, influencing the range of influence of a single training example. The \( \epsilon \)-insensitive tube, which defines the margin of tolerance within which no penalty is given to errors, is configured at 0.1. 
    
    \textbf{Random Forest:} The Random Forest model is instantiated with 100 trees, ensuring a good balance between performance and overfitting. The maximum depth of each tree is limited to 10 to prevent overfitting by controlling the complexity of the decision paths. 
    
    \textbf{XGBoost:} The XGBoost configuration involves setting the maximum depth of each tree to 5, which limits the model complexity and helps in reducing overfitting. The learning rate is set at 0.1, which specifies the step size at each iteration and helps in controlling the speed of convergence. The objective function is specified as a regression with squared error as the loss to be minimized. The model undergoes training for 100 iterations, with each iteration incrementally improving the model based on the defined learning rate and depth.

    The results indicate that the Inception Network outperforms the other three benchmark models in predicting all CICY triple intersection numbers, especially in $d_3$ and $d_p$, which have more complex data structures and harder to predict. Our Inception Network model demonstrates significant advantages, achieving substantially lower MSE compared to the other models. This observation suggests that the Inception Network more effectively captures and models the complex data structures associated with CICY manifolds.
        
    Furthermore, in Figure \ref{fig:different methods for different radio of training data}, we compare the performance of our Inception Network model against traditional machine learning regressors like SVM regressor, Random Forest, and XGBoost, across varying training set sizes. Notably, our method demonstrates remarkable stability even with smaller training datasets, as evidenced by the consistently high accuracy and low MSE across the fractions of data used for training. When assessing accuracy, the Inception Network shows a steady near-plateau performance from a training fraction of approximately 0.5 onwards, achieving almost peak accuracy, which is substantially higher than the other methods which either increase gradually or plateau at lower levels. Similarly, in terms of MSE, our Inception Network not only starts with lower error rates at smaller data fractions but also maintains a consistent decrease in errors as more data is introduced. This consistency in lower MSE and higher accuracy underscores the robustness of the Inception Network. The results illustrate the advantages of advanced neural network architectures in learning complex patterns more effectively than conventional methods.
    
    We note that for all baseline regressors (SVM, Random Forest, XGBoost), we did not perform extensive hyperparameter tuning beyond standard package defaults or common settings. While this provides a fair benchmark comparison, their performance could likely be further improved with systematic hyperparameter optimization.

\section{Discussion}

In this paper, our aim is to harness deep learning, specifically the Inception Network model, to predict the triple intersection numbers of CICY manifolds, enhancing the application of such advanced computational techniques to geometry.

Our results show the Inception Network's superiority in  performance for the CICY triple intersection numbers $d_1$, $d_2$, $d_3$ and $d_p$ as defined in \eqref{d-def}. The significance of our study lies in the high accuracy and low MSE achieved, pointing to the robustness of the Inception Network in handling topology datasets. This performance surpasses traditional regression models. 

Furthermore, our model demonstrates an advancement over previous works focused on Hodge number prediction, suggesting a broader potential for deep learning applied to algebraic geometry. It establishes deep learning's efficacy in capturing complex mathematical structures, which can be generalized to other areas of AI for science. However, we acknowledge the model's limitations, such as in predicting $d_p$, which presents a challenge due to its more sophisticated construction. 

We recommend future research to further refine the model's predictive capabilities including utilizing the recent developed manifold fitting techniques \citep{yao2019manifold,yao2023manifold,yao2023pnas,yao2024pnas}, explore additional properties of CICY manifolds, and extend the application of deep learning models to a wider array of mathematical structures. 
For instance, a nice work \citep{coates2023machine} explored the quantum dimension of Fano varieties, it would be a very fruitful study to be able to predict more refined objects such as Gromow-Witten invariants of our CICY manifolds.

The ultimate goal of this study, as with all AI-driven mathematical discovery \citep{He:2019nzx,He:2024gnk}, is to uncover new, interpretable structures, either in the form of a conjectured formula or a new derivation/proof. A clear future work is to extract a conjectural formula amongst $d_{1,2,3,p}$ in an automated way, that would pass the Automaticity, Interpretability, and Non-triviality of the Birch test of AI-driven findings \citep{he2024can}.

\bibliographystyle{plainnat}
\bibliography{references}

\end{document}